\documentclass[12pt,leqno,twoside]{amsart}
\usepackage{amssymb}
\usepackage{latexsym}

\usepackage{epsfig}
\usepackage{changebar}
\usepackage{srcltx}

\newtheorem{theorem}{Theorem}[section]
\newtheorem{corollary}[theorem]{Corollary}
\newtheorem{lemma}[theorem]{Lemma}
\newtheorem{proposition}[theorem]{Proposition}
\theoremstyle{remark}
\newtheorem{remark}[theorem]{\sc Remark}
\theoremstyle{remark}

\theoremstyle{definition}
\newtheorem{definition}[theorem]{Definition}
\theoremstyle{remark}
\newtheorem{example}[theorem]{\sc Example}

\theoremstyle{remark}

\theoremstyle{remark}

\renewcommand{\Box}{\square}    %\diamond

\newcommand{\cal}{\mathcal}

\renewcommand{\int}{{\rm{int}}}

\newcommand{\Sing}{{\rm{Sing\hspace{2pt}}}}

\newcommand{\mult}{{\rm{mult}}}

\newcommand{\cl}{{\rm{closure}}}

\newcommand{\ity}{{\infty}}

\newcommand{\e}{\varepsilon}
\newcommand{\m}{\setminus}
\newcommand{\fin}{\hspace*{\fill}$\Box$\vspace*{2mm}}

\newcommand{\cG}{{\cal G}}
\newcommand{\cH}{{\cal H}}

\newcommand{\cO}{{\cal O}}

\newcommand{\cN}{{\cal N}}

\newcommand{\bC}{{\mathbb C}}

\newcommand{\bP}{{\mathbb P}}

\newcommand{\bV}{{\mathbb V}}
\begin{document}

\title[Betti bounds of polynomials]
 {Betti bounds of polynomials}

\author{\sc Dirk Siersma}  

\address{Institute of Mathematics, Utrecht University, PO
Box 80010, \ 3508 TA Utrecht
 The Netherlands.}

\email{D.Siersma@uu.nl}

\author{\sc Mihai Tib\u ar}

\address{Math\' ematiques, UMR 8524 CNRS,
Universit\'e Lille 1, \  59655 Villeneuve d'Ascq, France.}

\email{tibar@math.univ-lille1.fr}

\thanks{The authors acknowledge the support of the Mathematisches Forschungsinstitut Oberwolfach under the frame of the Research-in-Pairs program.}



\date{January 17, 2011}

\dedicatory{Dedicated to the memory of Vladimir Arnol'd}

%\commby{}

%%% ----------------------------------------------------------------------

\begin{abstract}
  We initiate a classification of polynomials $f : \bC^n \to \bC$ of degree $d$ having
the top Betti number of the general fibre close to the maximum.
We find a range in which the polynomial must have isolated singularities and 
another range where it may have at most a line singularity of Morse transversal type, besides controlled singularities at infinity. Our method  uses deformations into particular pencils with non-isolated singularities.

\end{abstract}

%%% ---------------------------------------------------
\maketitle
%%% ---------------------------------------------------

\setcounter{section}{0}

\section{Introduction and results}\label{s:intro}

Let $f: \bC^n \to \bC$ be a polynomial function of degree $d\ge 2$, where $n\ge 2$. It is well-known that $f$ is a locally trivial fibration over $\bC$ outside a finite number of \textit{atypical values}, \cite{Th, Br}. It's \textit{general fibre} $G$ is a 
Stein manifold and therefore homotopy equivalent to a CW-complex of dimension $n-1$, by \cite{Ka, Ha}.
  We shall call \textit{top Betti number of $f$} and denote $b_{n-1}(f) := b_{n-1}(G)$ the $(n-1)$th Betti number of the general fibre. While this number is clearly  bounded in terms of $n$ and $d$, our aim is to find what are the special properties of $f$ which make  $b_{n-1}(f)$ approach the maximum $(d-1)^n$. Let us call \textit{top Betti defect of $f$} the  difference  $\Delta_{n-1}(f):= (d-1)^n - b_{n-1}(f)$.  

We prove in this paper that if $\Delta_{n-1}(f)$ is small enough, then the polynomial has special types of singularities. To do that, we study deformations of polynomials from a new viewpoint, namely by following the general fibre in the deformation process and drawing consequences on the ``quantity'' of singularities measured by the top Betti number.

Deformations of polynomial functions with given degree have been studied from the point of view of topological equivalence. In the 1960s Thom conjectured, and then Fukuda \cite{Fu} proved, that there are finitely many equivalence classes. In  \cite{BT} one employs previous results on deformation of germs to extend the L\^e-Ramanujam-Timourian criterion for the topological equivalence in families of germs to families of polynomials by using the constancy of Milnor numbers including those at infinity. 
The deformation theory of projective hypersurfaces, a classical topic with more recent contributions by Shustin and Tyomkin \cite{ShT}, is also related to the deformation of polynomials. In \cite{SS} one uses deformations to classify polynomials of 2 variables and degree 4, with respect to their singularities at infinity.

In order to state the precise results, we need a few definitions.
We consider $\bP^n$ as the standard compactification of $\bC^n$ for some fixed affine coordinates. We use the following notations:
  $X_t = f^{-1}(t)$ is some fibre of $f$ and $\overline{X_t} \cap H^\ity$ is the intersection of its closure in $\bP^n$ with the hyperplane at infinity $H^\ity := \bP^{n-1}$. We denote by $f_d$ and $f_{d-1}$ the degree $d$ and $d-1$ parts of $f$, respectively. Let  $\Sigma_f^\ity  := \Sing f_{d} = \{ [x] \in \bP^{n-1} \mid \frac{\partial f_d}{\partial x_i}(x) =0, i= 1,\ldots , n \}$. Geometrically, $\Sigma_f^\ity$ is the set of singular points of the restriction $\overline{X_t}\cap H^\ity$ and this set does not depend on the value $t\in \bC$. We call such points \textit{tangencies at infinity}. 
%\begin{definition}%\label{d:genericatinfinity}
An affine hypersurface in $\bC^n$ will be called \textit{general-at-infinity} if its projective closure is non-singular in the neighbourhood of the hyperplane at infinity $H^\ity$ and intersects it transversely. 
The polynomial $f$ will be called {\em general-at-infinity} (or \textit{of $\cG$-type}) iff all its fibres are general-at-infinity. Since this is equivalent to $\Sigma_f^\ity = \emptyset$, it means that $f$ is  general-at-infinity iff some fibre of $f$ is so.

%We show that in the next range of width $d-1$ the polynomial may have only special types of %non-isolated singularities, as follows:
 %%%%%%%%%%%%%%%%%%%%%%%%
\begin{theorem}\label{t:bettimax} 
Let $f : \bC^n \to \bC$ be some polynomial of degree $d \ge 2$. Then:
\begin{enumerate}
 \item  $\Delta_{n-1}(f) \ge 0$ 
and the equality holds if and only if $f$ is general-at-infinity.
\item If $0< \Delta_{n-1}(f) \le d-1$ then $\dim \Sing f \le 0$ and $\dim \Sigma_f^\ity \le 0$.
\item If  $d \le \Delta_{n-1}(f) < 2d-2$ for $d\ge 3$, then $\dim \Sigma^\ity_f \le 0$
and either $\dim \Sing f\le 0$ or $\Sing f$ is a line with generic Morse transversal type and transverse to the hyperplane at infinity.
\end{enumerate}
\end{theorem}
%%%%%%%%%%%%%%%%%%%%%%%

%\label{t:nonisolated2}

Let us remark that the above bounds depend only on the degree and not on the number of variables.
In our proof we have to start with some polynomial $f$ and reduce the dimensions of $\Sing f$ and of $\Sigma_f^\ity$ to 1, by using special deformations. Then we analyse one-dimensional singularities.
It appears for instance that the complementary case ``$\dim \Sing f\le 0$ and  $\dim \Sigma^\ity_f =1$''
does not occur in the range of Theorem \ref{t:bettimax}(c). To investigate this setting, we study in  \S \ref{s:onedimatinfinity} the Euler characteristic of a projective hypersurface with one-dimensional singular locus by developing a pencil technique which goes beyond the usual Lefschetz method since in our setting we deal with non-isolated singularities. This section is therefore of an independent interest. From its main result Theorem \ref{t:euler_sum}  we derive in \S \ref{s:line} that positive dimensional singularities in the hyperplane at infinity 
 produce a deeper Betti defect, namely:

\begin{theorem}\label{t:lineinfinity}
 If  $f$ has at most isolated affine singularities and $\Sigma^\ity_f$ is a reduced projective line with generic Morse transversal type, then
\[  \Delta_{n-1}(f) \ge 2(n-1)(d-2) +1. \]
\end{theorem}

These results set the bases of a classification of polynomial functions with top Betti number close to the maximum in terms of their singularities, and more particularly of those appearing at infinity. Our study is based on algebro-topological deformation tools and on the refined interplay between isolated and non-isolated singularities via polar curve techniques like explored in \cite{Le, Yo, Si1, Si2, Si3} etc.  

\smallskip

In case of polynomials with isolated singularities, namely $\dim \Sing f \le 0$ and $\dim \Sigma_f^\ity \le 0$, one may list the possible
combinations of critical points at infinity which yield generic fibres with top
Betti number close to the maximum. We give below the list for any $n\ge 2$, up to the top Betti defect equal to 3.   The notation $\langle \cdot |\cdot \rangle$ means the singularity type of the boundary pair $(\overline{G}, \overline{G} \cap H^\ity)$
at some point at infinity (see \S \ref{s:prel}), where $G$ is some general fibre of $f$. The additive notation is used if multiple
special points occur at infinity. On the second column there appear the singularity types at one, two or three points at infinity, where $A_0$ stands for ``nonsingular''. The third column contains the Arnol'd type fraction notation, cf \cite{Ar}. According to Theorem \ref{t:bettimax}(b), this gives all the types of general fibres of polynomials with small top Betti defect. The realisation of each type is an easy exercise; note that for some classes we need high enough degrees.

\[
\begin{array}{|l|l|l|}
%\hline
%\multicolumn{3}{|c|}{\mbox{F polynomials}}\\
\hline
\mbox{$\Delta_{n-1}(f)$} & \mbox{boundary type} & \mbox{Arnol'd type} \\
\hline
0     & \langle A_0|A_0\rangle &   A_0 \\
\hline
1 & \langle A_0|A_1\rangle &   A_1 \\
 \hline
2 & \langle A_0|A_2\rangle & A_2\\
            &  2\langle A_0|A_1\rangle & 2 A_1\\
            &  \langle A_1|A_1\rangle & B_2\\
 \hline
 3 & \langle A_0|A_3\rangle & A_3\\
            & \langle A_0|A_2\rangle + \langle A_0|A_1\rangle & A_2 + A_1\\
            &  3\langle A_0|A_1\rangle & 3 A_1\\
            & \langle A_1|A_2\rangle & C_3 \\
            & \langle A_1|A_1\rangle + \langle A_0|A_1\rangle & B_2 + A_1\\
            & \langle A_2|A_1\rangle & B_3 \\
\hline	    
\end{array}
\]

\medskip

\noindent 
Let us remark that \cite{SS} contains a certain algebraic classification for $n=2$ and $d=4$, together with the computation of the Milnor numbers of the singularities at infinity, continuing the classification by C.T.C. Wall in degree 3 and two variables \cite{Wa}.

\smallskip

We close this introduction by a few examples, some of them showing that the bounds in our theorems are sharp.
\begin{example}\label{ex1}
 $f=x+ x^2 y$. We have $n=2$, $d=3$, $\Sing f = \emptyset$ and $\Sigma^\ity_f =[0;1]\in \bP^1$. The computation yields $b_1(f) =1$ and therefore $\Delta_{1}(f) = 4-1 = 3$ which corresponds to the situation in Theorem \ref{t:bettimax}(c).
\end{example}
\begin{example}\label{ex2}
 $f=x^2 y$. Here $n=2$, $d=3$, $\Sing f = \{ x=0\}$ is a line with transversal type $A_1$ and $\Sigma^\ity_f =[0;1]\in \bP^1$. By computation we have $b_1(f) =1$ and so $\Delta_{1}(f) = 4-1 = 3$. This corresponds  to the situation in Theorem \ref{t:bettimax}(c) and also shows that the estimation is sharp. In full generality, for any $n\ge 2$ and $d>2$, let
$f = (a_1 z_1^2 + \cdots + a_{n-1} z_{n-1}^2)x^{d-2} +  c_1 z_1^d + \cdots + c_{n-1} z_{n-1}^d$.
This is a homogeneous polynomial and has a line singularity $L = \{ z_1 = \ldots = z_{n-1} =0\}$. For general coefficients $a_i$, $c_i$, this polynomial has no other singular point and the line $L$ is transversal to $H^\ity$ and has Morse generic transversal singularity type. The projective hypersurface $\{ f_d =0 \}$ has a single singular point at $p:= [1; 0; \cdots ; 0]$.
By a local computation, the singularity type of $\overline{X_t}$ at $p$ is $A_{d-1}$ for $t\not=0$.
Then we may apply formula (\ref{eq:F}) of Proposition \ref{p:defect} since by Remark \ref{r:defect}
this works in our situation too. Indeed, the fibre $X_t$ of $f$ has reduced homology concentrated in dimension $n-1$ since it is diffeomorphic to the Milnor fibre of $f$ at the origin (by the homogeneity) and since by \cite{Si1} this line singularity Milnor fibre is homotopy equivalent to a bouquet of spheres of dimension $n-1$. So by (\ref{eq:F}) we get $\Delta_{n-1}(f) = d-1 + 1 = d$, which also shows that the lower bound in Theorem \ref{t:bettimax}(c) is sharp.
\end{example}

 \begin{example}\label{ex4}
 Let $f= z^d + z^2x^{d-2} + z^2y^{d-2} + xy(x^{d-3} - y^{d-3})$, where $n=3$ and $d\ge 3$. Then $f$ has isolated affine singularities and $\Sigma^\ity_f$ is a reduced projective line with generic transversal type $A_1$. Let us compute $\Delta_{2}(f)$ using Proposition \ref{p:defect} and its notations.

 \noindent
Let $X_t$ denote a general fibre. Then $H^\ity \cap \Sing \overline{X_t} = \{ z=0, xy (x^{d-3} - y^{d-3})=0\}$.
By computing in coordinate charts it follows that the set $\Sing \overline{X_t}$ consists of $d-1$ Morse singularities. According to formula (\ref{eq:B}), this contributes with $d-1$ to the top Betti defect. 
Next we compute $\Delta \chi^\ity$ from the same formula. We have $f_d = z^2(z^{d-2} + x^{d-2}+ y^{d-2})$ and the reduced hypersurface $\{f_d =0\} \subset \bP^2$ has degree $d-1$ and precisely $d-2$ singular points which are Morse. This implies the equality $\chi(\{f_d =0\}) = \chi^{2,d-1} + d-2$. Since $\chi^{2,d-1} = -(d-1)^2 + 3(d-1)$, we get $\chi(\{f_d =0\}) = -d^2 + 6d-6$, which yields $\Delta \chi^\ity = \chi^{2,d} -\chi(\{f_d =0\}) = -3d + 6$. By formula  (\ref{eq:B}) we then get $\Delta_{2}(f) = d-1 + 3d - 6 = 4d -7$. 
 This corresponds to the equality in Theorem \ref{t:lineinfinity}, showing that the bound is sharp for $n=3$. See also Remark \ref{r:ex}.
\end{example}

%\smallskip

%%%%%%%%%%%%%%%%%%%%%%%%%%%%%%%%%%%%%%%%%%%%%%%%%%%%%%%%%%%%%%%

%%%%%%%%%%%%%%%%%%%%%%%%%%%%%%%%%%%%%%%%%%%%%%%%%%
%%%%%%%%%%%%%%%%%%%%%%%%%%%%%%%%%%%%%%%%%%%%%%
%%%%%%%%%%%%%%%%%%%%%%%%%%%%%%%%%%%%%%%%%%%%%%%%%%%%

\section{Preliminaries on deformations of polynomials}\label{s:prel}

By \textit{one-parameter deformation of $f$} we mean a holomorphic map $P : \bC^n \times \bC \to \bC$, where  $P_s:= P(\cdot, s)$ is a polynomial of degree $d$ for any  $s\in \bC$, and such that $P_0 = f$.  We shall work with germs at $s=0$ of such families of polynomials. Let $G_s$ denote the general fibre of $P_s$. We start from the following result about the behaviour of the general fibre in such a deformation.

\begin{proposition}\label{t:spec}  \cite{ST-ExCh}
For $s\not= 0$ close enough to $0$,
 the general fibre $G_0$ of $P_0$ can be naturally embedded in the 
general fibre $G_s$ of $P_s$ such that the embedding $G_0\subset G_s$ induces an injective map $H_{n-1}(G_0)
 \hookrightarrow  H_{n-1}(G_s)$.
\fin
\end{proposition}
One derives from this the following \textit{semi-continuity principle} for the top Betti number,  a cornerstone of our study:
\begin{equation}\label{c:globsemi}
\Delta_{n-1}(P_s) \le \Delta_{n-1}(P_0), \mbox{ for } s\not= 0 \mbox{ close enough to } 0.
\end{equation}

  Another key fact, that we shall now prove, is that any polynomial is deformable into a general-at-infinity polynomial (defined in \S \ref{s:intro}).

It is well-known that the Euler characteristic of a non-singular and general-at-infinity hypersurface of degree $d$ in $\bC^n$ is equal to $1 + (-1)^{n-1} (d-1)^{n}$.
By its definition, a $\cG$-type polynomial may have special fibres with at most isolated singularities.
An example of a general-at-infinity polynomial of degree $d$ is $x_1^d + \cdots + x_{n}^d$.

\begin{proposition}\label{p:general}
Any polynomial can be deformed into a general-at-infinity polynomial of the same degree. 
 More precisely, let $h_d$ be some general-at-infinity polynomial of degree $d$. Then the deformation $f_\e:= f + \e h_d$ transforms any given polynomial $f$ of degree $d$ into a general-at-infinity polynomial $f_\e$, for any $\e \not=0$ close enough to $0$.
\end{proposition}
\begin{proof}

Let us first remark that one may deform any hypersurface $X\subset \bC^n$
 in a constant degree family $\{ X_s\}_{s\in \delta}$ such that  $X=X_0$ and that $X_s$ is
general-at-infinity and nonsingular for $s\not= 0$ in a small enough disk $\delta$ centered at $0\in \bC$, as follows. Let $X_s = \{ (1-s)f + s(h_d -1) = 0\}$.  The family $\{X_s=0\}_{s\in [0,1]}$ has a finite number of special members since it contains the general-at-infinity non-singular hypersurface $\{ h_d = 1$\} and since transversality and non-singularity are open properties. It follows that only a finite number of the hypersurfaces of this pencil are not general-at-infinity or singular.
Then the family $\{ X_s\}_{s\in \delta}$  has the desired property for a small enough disk $\delta$ centered at $0$. 

 Next, consider the deformation $f_\e = f + \e h_d$ of $f$. Taking $s= \frac{\e}{\e-1}$ in the above pencil we deduce that the fibre $f_\e = \e$ is general-at-infinity. Since the generic fibre of a polynomial is well-defined, we deduce that $f_\e$ is a general-at-infinity polynomial for any small enough $\e \not= 0$.
\end{proof}

%%%%%%%%%%%%%%%%%%%%%%%%%%%%%%%%%%%%%%%%%%%%%%%%%
%%%%%%%%%%%%%%%%%%%%%%%%%%%%%%

%Let us recall from \cite{Di, ST-gendefo, Ti} some facts about the topology of the general fibre of a %polynomial $f$ in terms of the boundary singularities of its projective compactification. 

\subsection*{Boundary singularities}
We assume in the remainder of this paper that $f : \bC^n \to \bC$ is a polynomial function of degree $d$. Let $Y_t = f^{-1}(t)$.
At some point $p \in H^\ity$ we consider the \textit{boundary germ} $(\overline{Y_t},\overline{Y_t} \cap H^\ity)_p$ which is actually a family of germs depending on $t\in \bC$. This family has a constant singularity type outside finitely many values of $t$. The singularity theory of germs of boundary singularities with respect to a hyperplane has been studied by Arnol'd \cite{Ar-bd}. He introduced the concept of isolated boundary singularity and gave a list of simple singularities. Some of the types occur for instance in our table in \S \ref{s:intro}.
 
\begin{definition} 
One says that the boundary pair $(\overline{Y_t},\overline{Y_t} \cap H^\ity)_p$ has an \textit{isolated singularity} if  both $\overline{Y_t}$ and $\overline{Y_t} \cap H^\ity$ have (at most) isolated singularities at $p$.
\end{definition}

We have the equivalence: the boundary pair $(\overline{Y_t},\overline{Y_t} \cap H^\ity)$ has isolated singularities if and only if $Y_t$ has isolated singularities and $\dim \Sigma_f^\ity \le 0$. Then  $\Sigma_f^\ity \cap \{f_{d-1}= 0\}$ is the subset of points of $H^\ity$ where $\overline{Y_t}$ is singular, and this does not depend on the value $t\in \bC$.

\begin{proposition}\label{p:defect} 
Let $f$ be a polynomial of degree $d$ and isolated singularities, having general fibre $X_0$ and satisfying $\dim \Sigma_f^\ity \cap \{f_{d-1}= 0\} \le 0$. Then:
\begin{equation}\label{eq:B}
\Delta_{n-1}(f) = \sum_{p\in\Sigma_f^\ity\cap \{f_{d-1}= 0\} }\mu_p(\overline{X_0}) + (-1)^n \Delta \chi^\ity
\end{equation}
where $\Delta\chi^{\infty} := \chi ^{n-1,d} - \chi (\{ f_d = 0 \})$ and
$\chi ^{n-1,d}  =
  n  - \frac{1}{d} \{ 1 + (-1)^{n-1} (d-1)^{n} \}$ denotes
the Euler characteristic of the smooth hypersurface $V^{n-1,d}_{gen}$ of degree $d$
in $\bP^{n-1}$.

In particular, if  $\dim \Sigma_f^\ity \le 0$ then:
\begin{equation}\label{eq:F}
 \Delta_{n-1}(f) = \sum_{p\in \Sigma_f^\ity } [ \mu_p(\overline{X_0})+ \mu_p(\overline{X_0} \cap H^\ity)].
\end{equation}
\end{proposition}

\begin{proof}
Formula (\ref{eq:B}) was stated without proof in \cite[(6.1 and (6.2)]{ST-ExCh}, \cite[\S 4]{Ti}.  We consider the germ of a deformation of $X_0$ in a constant degree family $\{ X_s\}_{s\in \delta}$ such that $X_s$ is non-singular and
general-at-infinity for $s\not= 0$; this exists according to Proposition \ref{p:general}.
For any $s$, we have the equality of Euler characteristics $\chi(X_s) = \chi(\overline{X_s}) - \chi(\overline{X_s} \cap H^\ity)$ and taking the difference we get, for $s\not=0$:
\[
 \chi(X_0) - \chi(X_s) = \chi(\overline{X_0}) -  \chi(\overline{X_s}) - \chi(\overline{X_0} \cap H^\ity) + \chi(\overline{X_s} \cap H^\ity)
\]
 Our family $\{ \overline{X_s}\}_{s\in \delta}$ is in particular a smoothing of the hypersurface $\overline{X_0}$ with isolated singularities and thus the jump of Euler characteristics is the sum of the Milnor numbers of the singularities of $\overline{X_0}$:
\[
 \chi(\overline{X_0}) -  \chi(\overline{X_s}) = (-1)^n \sum_{p\in \Sigma_f^\ity \cap \{f_{d-1}= 0\}} \mu_p(\overline{X_0}).
\]
In case  $\dim \Sigma_f^\ity \le 0$, the family $\overline{X_s} \cap H^\ity$ is a smoothing of 
$\overline{X_0} \cap H^\ity$ and we get the following similar equality: $\chi(\overline{X_0}\cap H^\ity) -  \chi(\overline{X_s}\cap H^\ity) = (-1)^{n-1} \sum_{p\in \Sigma_f^\ity} \mu_p(\overline{X_0}\cap H^\ity)$.

Up to this point the Euler characteristic computation runs like in \cite[ch.1, Prop. 4.6]{Di}. Next we need that the general fibre
 has concentrated homology in the top dimension, which is true since the polynomial $f$ has isolated singularities in the sense of \cite{Pa} or \cite{ST} and by \textit{loc.cit.}, the general fibre has the homotopy type of a bouquet of spheres $S^{n-1}$ and therefore its reduced homology is concentrated in dimension $n-1$. This implies that $\chi(X_0) = 1 - (-1)^nb_{n-1}(X_0)$.  We also have that $\chi(X_s) = 1- (-1)^n(d-1)^n$ since $X_s$ is general-at-infinity and non-singular.
Collecting these information we get formulas (\ref{eq:B}) and (\ref{eq:F}). 
\end{proof}

\begin{remark}\label{r:defect}
 The above proof shows that the hypothesis ``$\dim \Sing f \le 0$'' of Proposition \ref{p:defect} may be replaced by ``the general fibre of $f$ has the reduced homology concentrated in dimension $n-1$''.
\end{remark}

%%%%%%%%%%% hypersurfaces

%%%%%%%%%%%%%%%%%%%%%%%%%%%%%%%%%%%%%%%%%%%%%%%%%%%%%%%%%%%%%%%%%

\section{Polynomials with non-isolated singularities}

Let $f : \bC^ n \to \bC$ be a polynomial of degree $d\ge 2$ with singular loci of dimension $\ge 2$, more precisely $\dim \Sing f \ge 2$ or $\dim \Sigma_f^\ity \ge 2$. If one deforms $f$ directly to general-at-infinity polynomials, then it appears that comparing the general fibres becomes a difficult task. A better strategy would be to deform $f$ in two steps and use the semi-continuity principle (\ref{c:globsemi}) according to the following program:
(a). deform such that the dimension of the singularity locus decreases to one,  and then (b). compare the new polynomial to another deformation of it into a polynomial satisfying the hypothesis of Proposition \ref{p:defect} or directly use results from the theory of one-dimensional singularities. 
The reason is that one-dimensional singularities and their deformations are quite well understood, due to the work of L\^e \cite{Le}, Yomdin \cite{Yo} and the detailed study by Siersma \cite{Si1} and his school, see e.g. the survey \cite{Si-Cam}. 

Let $l: \bC^ n \to \bC$ be a linear function.
We denote by:
\[ \Gamma (l,f) := \cl [\Sing (l,f)\setminus \Sing f ] \subset \bC^n\]
the \textit{polar locus of $f$ with respect to $l$}. One has the following Bertini type result, proved in \cite{Ti-compo, Ti-equi}, \cite[Thm. 7.1.2]{Ti}:

\begin{lemma}\label{t:bertini} 
There is a Zariski-open subset $\Omega_{f}$ of the dual projective space
 $\check\bP^{n-1}$ such
that, for any $l\in \Omega_f$, the polar locus $\Gamma(l,f)$  is a reduced curve or it is
empty. \fin
\end{lemma}
We may and shall also assume (by eventually restricting $\Omega_f$ to some open Zariski subset of it) that if $\dim \Sing f \ge 1$ then $\dim \Sing f \cap \{ l=0\} = \dim \Sing f -1$ for any $l\in \Omega_f$. We then say that $l$ is \textit{general with respect to $f$} whenever $l\in \Omega_f$. With these settings we may start our program.

\begin{lemma}\label{l:1}
 Let $l$ be general with respect to $f$ and to $f_d$. If $\dim \Sing f \ge 1$, or if $\dim \Sigma_f^\ity \ge 1$, then the deformation $f_\e = f + \e l^d$
reduces by one the dimension of the respective singular locus. If $f$ has the property that $\dim \Sing f \le 0$ or $\dim \Sigma^\ity \le 0$ then $f_\e$ preserves this property.
\end{lemma}
\begin{proof}
 After some linear change of coordinates, we may assume that $l=x_n$. Then $\Sing (l,f) = \{ \frac{\partial f}{\partial x_1}=0, \ldots , \frac{\partial f}{\partial x_{n-1}} =0 \}$. Let us first show that
$\dim \Sing f_\e =  \dim \Sing f -1$ whenever $\dim \Sing f \ge 1$.
We have the inclusions:
\[ \Sing f \cap \{ l=0\} \subset \Sing f_\e \subset (\Sing f\cap \{ l=0\}) \cup \Gamma (l,f)
\]
since if $x\in \Sing f_\e \cap \Sing f$ then $x\in \{ l=0\}$.
The algebraic function $\frac{\partial f}{\partial x_n} + \e d l^{d-1}$ cannot be identically zero on some component $\Gamma_i \subset  \Gamma(l,f)$ for more than one nonzero value 
of $\e$. These facts show that, for $\e\not= 0$ in some small enough neighbourhood of $0$ in $\bC$,  $\Sing f_\e$ is the union of $\Sing f\cap \{ l=0\}$ and a finite set of points.

Next, let us show what is the effect of this deformation on $\Sigma_f^\ity$.  Taking $f_d$ instead of $f$ in the above proof, the genericity of $l=x_n$ implies, via Lemma \ref{t:bertini}, that $\Gamma(l,f_d)$ is either empty or a homogeneous algebraic set of dimension 1 in $\bC^n$ and hence a finite set of points in $\bP^{n-1}$. It is also clear that $\Sing (f_d + \e l^{d}) \cap \Sing f_d = \Sing f_d \cap \{ l=0\}$. These show that $\Sing (f_d + \e l^{d})$ is the union of $\Sing f_d \cap \{ l=0\}$ and a finite set of points (eventually empty) and therefore that $\dim \Sigma_{f_\e}^\ity = \dim\Sigma_f^\ity -1$, provided that the former is $\ge 1$.
\end{proof}

\begin{lemma}\label{l:2}
Let $l$ be general with respect to $f$ and consider the deformation $f_\e = f + \e l$. If $\dim \Sing f \ge 1$ then there exists a small disk centered at the origin $D\subset \bC$ such that $\dim \Sing f_\e \le 0$ and $\Sigma_{f_\e}^\ity = \Sigma_f^\ity$, for any $\e \in D^*$.
\end{lemma}
\begin{proof}
As in the above proof we may assume that $l=x_n$. For any $\e \not= 0$ we have $\Sing f \cap \Sing f_\e = \emptyset$ and the polynomial function $\frac{\partial f}{\partial x_n} + \e$ cannot be identically zero on some irreducible component $\Gamma_i \subset  \Gamma (l,f)$ for two different values of $\e$. The claim follows.
\end{proof}

With these preparations we may consider in the next statements the two cases of dimension one singular locus.

\begin{proposition}  \label{p:01}
If $f$ is a polynomial of degree $d$ with $\dim \Sigma^\ity_f = 1$ and   $\dim \Sing f \le 0$ then
 \[ \Delta_{n-1}(f) \ge d-1.\]
\end{proposition}
\begin{proof}
 We consider the deformation $f_\e = f + \e l^d$ for general $l$ as in Lemma \ref{l:1}. It then follows that $\dim \Sing f_\e \le 0$ and that $\Sing f_{\e, d} \subset \bP^{n-1}$ is the union of $\Sigma^\ity_f \cap \{ l=0\}$ and eventually some finite set of points. We may assume that the hyperplane $\{ l=0\}$ slices $\Sigma^\ity_f$ at regular points only; this property is generic too. Let then $\Sigma^\ity_f = \cup_r \Sigma_r$ be the decomposition into irreducible components and let $p\in \Sigma_r \cap \{ l=0\}$ for some $r$. In order to compute  the top Betti defect $\Delta_{n-1}(f_\e) = (d-1)^n - b_{n-1}(X_\e)$ for some general fibre $X_\e$ of $f_\e$, we may use formula (\ref{eq:F}) in which one of the ingredients
is $\mu_p (\overline{X_\e} \cap H^\ity)$ and observe that this is equal to $\mu_p (f_{\e, d})$.

 We denote by $\mu^\pitchfork_r$ the Milnor number of the transverse singularity of $\Sigma_r$.
By the local L\^ e \textit{attaching formula}, see \cite{Le}, \cite{Yo} and  \cite{Si3}, we have:
\begin{equation}\label{eq:le}
  \mu_p (f_{\e, d}) + \mu_p ({(f_{\e, d })}_{|l=0}) = \mult_p (\hat\Gamma_p(l, f_{\e, d }), \{ f_{\e, d } =0\} ) + \mult_p (\Sigma_r, \{ f_{\e, d } =0\})  \mu^\pitchfork_r 
\end{equation}
where $\hat\Gamma_p(l, f_{\e, d })$ denotes the union of the components of the germ at $p$ of the polar curve of the map $\psi := (l, f_{\e, d }) : (\bC^n, 0) \to (\bC^2, 0)$ other than the singular locus $\Sigma_r$.

\begin{figure}[hbtp]
\begin{center}
%\epsfxsize=5cm
%\leavevmode
\resizebox{!}{1.25in}{\includegraphics{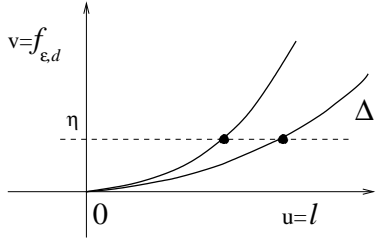}}
\end{center}
\caption{{\em Polar multiplicities
 }}
\label{f:1}
\end{figure}

 By the regularity of $p$, it follows that $\mu_p ({(f_{\e, d })}_{|l=0}) = \mu^\pitchfork_r$.
In local coordinates at the regular point $p$ the germ of the singular locus $\Sigma_r$ is a line
and the restriction to $\Sigma_r$ of the map $\psi$ is one-to-one. The germ at $\psi(p)$ of the image $\Delta := \psi(\Sigma_r)$
is parametrised by $(l, \e l^d)$ since $f_{\e, d} = f_d + \e l^d$ and $\Sigma_r\subset \{ f_d =0\}$.
This multiplicity is represented in Figure \ref{f:1} by the number of intersection points of $\{ f_{\e, d } = \eta\}$ with the curve $\Delta$.  Therefore $\mult_p (\Sigma_r, \{ f_{\e, d } =0\}) = \mult_{\psi(p)}(\Delta, \{ v = 0\}) = d$, where $(u,v)$ are the coordinates of the target $(\bC^2, 0)$. Then formula (\ref{eq:le})  becomes:

\begin{equation}\label{eq:le2}
  \mu_p (f_{\e, d})  =  (d-1) \mu^\pitchfork_r  + \mult_p (\hat\Gamma_p(l, f_{\e, d }), \{ f_{\e, d } =0\} ) .
\end{equation}
We next need to sum up over all the points $p\in \Sigma^\ity_f \cap \{ l=0\}$. The number of points of $\Sigma_r \cap \{ l=0\}$ is equal to the degree $d_r := \deg \Sigma_r$ and we get:
\begin{equation}\label{eq:le3}
\sum_{p\in \Sigma^\ity_f\cap \{ l=0\}}  \mu_p (f_{\e, d})  =  (d-1) \sum_r d_r \mu^\pitchfork_r + \sum_r \sum_{p\in \Sigma_r\cap \{ l=0\}} \mult_p (\hat\Gamma_p(l, f_{\e, d }), \{ f_{\e, d } =0\} ) ,
\end{equation}
hence 
\begin{equation}\label{eq:le4}
 \sum_{p\in \Sigma^\ity_f\cap \{ l=0\}} \mu_p (f_{\e, d}) \ge  (d-1) \sum_r d_r \mu^\pitchfork_r 
\end{equation}
with equality if and only if $\hat\Gamma_p(l, f_{\e, d }) =\emptyset$ for all $p\in \Sigma^\ity_f\cap \{ l=0\}$.
We finally get from formulas (\ref{eq:F}) and (\ref{eq:le4}):
\begin{equation}\label{eq:le5} 
 \Delta_{n-1}(f_\e) = (d-1)^n  - b_{n-1}(f_{\e, d})  \ge  \sum_{p\in \Sigma^\ity_f \cap \{ l=0\}}  \mu_p(\overline{X_\e} \cap H^\ity)\ge  (d-1) \sum_r d_r \mu^\pitchfork_r \ge d-1,
\end{equation}
where $X_\e$ denotes the general fibre of $f_\e$.
The first inequality becomes an equality if and only if $\overline{X_\e}$ has no singularities in the neighbourhood of $H^\ity$. The last one becomes an equality if and only if $r=1$ and $d_1 = 1$.

Our claim follows since we have $\Delta_{n-1}(f) \ge \Delta_{n-1}(f_\e)$ by the semi-continuity principle (\ref{c:globsemi}).
\end{proof}
 
\begin{proposition}  \label{p:10}
If $f$ is a polynomial of degree $d$ with $\dim \Sing f = 1$ and $\dim \Sigma^\ity_f = 0$   then
\[ \Delta_{n-1}(f) \ge d-1.\]
\end{proposition}
\begin{proof}
 We may assume, by eventually adding up some appropriate constant, that the fibre $\{ f=0\}$ contains some 1-dimensional component of the singular locus $\Sing f$. If we consider the deformation $f_\e = f + \e l$ for general $l$ then, by Lemma \ref{l:2}, $f_\e$ has isolated singularities and $\Sigma^\ity_{f_\e} = \Sigma^\ity_{f}$  for small enough $\e \not=0$. We shall therefore use formula (\ref{eq:F}) in order to compute the top Betti defect $\Delta_{n-1}(f_\e)$, so let us take the homogenisation of degree $d$ by the variable $z$:
\[
 \{ \tilde f_\e = \tilde f + \e l z^{d-1}= tz^d\} = \overline{f_\e^{-1}(t)} \subset \bP^n.
\]
 As we did before, we may and shall assume that $l = x_n$.
Pick up some point $p\in \overline{f^{-1}(0)}\cap \overline{\Sing f}\cap H^\ity$ and note that $p\in \Sigma^\ity_{f} =  \Sigma^\ity_{f_\e}$. We may assume that $p\not\in \{ x_n =0\}$ and more specifically that $p=[0;0\cdots;0; 1]$. In the chart $x_n =1$ the function $f$ reads $\hat f$ and $\overline{f_\e^{-1}(t)}$ is defined by the equation $\hat f + (\e -tz) z^{d-1} =0$;  we consider its germ at the origin $0$ of the local coordinates. 
By taking some $(d-1)$th root $u= u(t,z)$ of the germ $\e -tz$ we change $z$ into $\hat z =uz$, which defines the same hyperplane at infinity.  The local equation at $p$ of the fibre $\overline{f_\e^{-1}(0)}$ becomes in our new coordinates:
\[
\hat f(x_1, \ldots, x_{n-1}, \frac{\hat z}{u}) + \hat z^{d-1} =0.
\]
 The polar locus $\Gamma_0(\hat z, \hat f + \hat z^{d-1})$ is a germ equal to the union $\Gamma_0(\hat f, \hat z) \cup \Sing_0 \hat f$.   Despite the fact that the coordinate $\hat z$ is maybe not generic in the sense of Lemma \ref{t:bertini}, $\Gamma_0( \hat z, \hat f)$ is a curve or empty by Pellikaan's result \cite[Prop 8.5]{Pe1}, since $\hat z$ has the property $\Sing \hat f \cap \{\hat z=0\} = \{0\}$. Let $\Sing_0 \hat f = \cup_j \Sigma_j$ be the decomposition into irreducible components of the germ and let us denote by $\mu^\pitchfork_j$ the Milnor number of the transversal singularity of $\Sigma_j$. 

 By formula (\ref{eq:F}) of Proposition \ref{p:defect}, the top Betti defect of $f_\e$ is expressed in terms of the boundary Milnor numbers. Each Milnor number may be expressed via the L\^e attaching formula in terms of polar multiplicity, see (\ref{eq:le}). Thus, by neglecting the contribution of the branches of $\Gamma_0( \hat z, \hat f)$ other than $\Sing_0 \hat f$, which is positive if and only if these branches are nonempty, we get:
\begin{equation}\label{eq:att1} 
 \mu_0(\hat f + \hat z^{d-1}) + \mu_0(\hat f_{| \hat z=0}) \ge  \sum_j \mult_0(\Sigma_j,  \{ \hat f + \hat z^{d-1} =0\}) \mu^\pitchfork_j.
\end{equation}

 Let $\Delta$ be the image of the singular locus  $\Sing_0 \hat f$ under the map $\psi := (\hat z, \hat f + \hat z^{d-1})$ and note that $\Delta$ is parametrised by $(t, t^{d-1})$.
The multiplicity $\mult_0 (\Sigma_j)$ is less or equal to $\mult_0 (\Sigma_j, \{z=0\})$ (with equality iff $\Sigma_j \pitchfork \{ z= 0\}$) and the later is the degree of the map $\psi_| : \Sigma_j \to \Delta$. We then have $ \mult_0(\Sigma_j,  \{ \hat f + \hat z^{d-1} =0\}) = (d-1) \deg_0 \psi_| \ge  (d-1) \mult_0( \Sigma_j)$ and thus we obtain 
 \begin{equation}\label{eq:att2} 
 \mu_0(\hat f + \hat z^{d-1}) + \mu_0(\hat f_{| \hat z=0}) \ge  (d-1)\sum_j \mult_0(\Sigma_j) \mu^\pitchfork_j.
\end{equation}
 
Since this inequality holds at any point $p\in \overline{f^{-1}(0)}\cap \overline{\Sing f}\cap H^\ity$, we may take the sum over all of those in order to bound from below the top Betti defect of $f_\e$. Since each term is a positive multiple of $d-1$, we get our claim.

\end{proof}

%%%%%%%%%%%%%%%%%%%%%%%%%%%%%%%%%%%%%%%%%%%%%%%%%%%%%%%%%%%%%%%%%%%%%

\section{Proof of Theorem \ref{t:bettimax}}\label{proof}
\subsection{Proof of part (b)} \ 
 Let $(\alpha, \beta) := (\dim \Sing f, \dim \Sigma_f^\ity)$. 
Let us observe that $\beta\ge \alpha -1$ since if some fibre of $f$ has a singular locus of dimension $\alpha' >0$ then it intersects the hypersurface at infinity $\{ f_d = 0\}$ and this intersection, which has dimension $\alpha'-1$, is included in the singular set of $\{ f_d = 0\}$ which is $\Sigma_f^\ity$ by definition.

  If $f$ has non-isolated singularities, i.e. $(\alpha,\beta) > (0,0)$ then, by applying Lemmas \ref{l:1} and \ref{l:2}, we deform $f$ until we arrive at some polynomial $f'$ which is in one of the following two terminal cases: $(\alpha, \beta) = (1,0)$ or $(0,1)$. More precisely the strategy is as follows. If $\beta >0$ then we first apply Lemma \ref{l:2} to get $(0,\beta)$, next apply Lemma \ref{l:1} repeatedly a number of $(\beta-1)$ times to reach the dimensions $(0,1)$.  In case $\beta=0$, the above observation shows that $\alpha\le 1$, hence we can only have $(\alpha,\beta) = (1,0)$.

 By the semi-continuity principle (\ref{c:globsemi}), we get  $\Delta_{n-1}(f)\ge \Delta_{n-1}(f')$. Then one applies either Proposition \ref{p:10} or Proposition \ref{p:01} to conclude in both cases that $\Delta_{n-1}(f') \ge d-1$. This shows that if $\Delta_{n-1}(f) \le d-2$ then $(\alpha,\beta) = (0,0)$.  Since our statement claims the inequality $\le d-1$, we need an extra argument; this will be given in \S \ref{ss:end1} bellow and will complete the proof.

%%%%%%%%%%%%%%%%%%%%%%%%%%%%%%%%%%%%%%%%%%%%%%%%%%%%%%%%%%%%%%%%%%%%%

\subsection{Proof of Part (c)}
As told above,  if $f$ has non-isolated singularities then we may deform $f$ until we may apply either Proposition \ref{p:01} or  Proposition \ref{p:10}.  In the first case we arrive at the 
the inequality (\ref{eq:le5}) which shows that, if for some $r$ we have either $r>1$, or $d_r >1$, or  $\mu^\pitchfork_r >1$,  then $\sum_r d_r \mu^\pitchfork_r \ge 2$ and therefore $(d-1)^n  - b_{n-1}(f_{\e, d})  \ge 2(d-1)$. In the second case we arrive at the 
the inequality (\ref{eq:att2}) which yields the same conclusion. 

Now if in the first case we have $r= d_r= \mu^\psi_r =1$ then we may use Theorem \ref{t:lineinfinity}, proved independently in \S \ref{s:line}, which tells that $\Delta_{n-1}(f) \ge 2(n-1)(d-2) +1$ and this is greater than $2d-3$ since $n\ge 3$. This completely eliminates the first case from our range $\Delta_{n-1}(f) < 2d-2$.

These arguments and the use of the range in Part (b) prove our statement. However we need to complete the proof of Part (b), and this will be done in the next.

\subsection{End of the proof of  Part (b)} \label{ss:end1}\ 
 Assume first that $(\alpha, \beta) = (1,0)$.
 The arguments presented in the above proof of Part (c) allow to deduce the following: if $\Delta_{n-1}(f) < 2d-2$ then $\Sing f$ is a line with transversal type $A_1$ intersecting $H^\ity$ transversely. Assume without loss of generality that $\Sing f\subset f^{-1}(0)$. Let then $\{p\} := \overline{\Sing f}\cap H^\ity$ and consider the local chart at this point, which we shall denote by $0$ in the following. 
Locally at $0$, the compactified fibres of $f$ have the equation $\hat f - t z^ d = 0$. Let us fix a general, not atypical value $t\not=0$. 
We may then apply the Yomdin-L\^e formula like in (\ref{eq:att1}) and (\ref{eq:att2}), but this time   with $d$ instead of $d-1$ and, by a completely similar computation, we get:
\begin{equation}\label{eq:att1bis}
  \mu_0(\hat f +  tz^{d}) + \mu_0(\hat f_{| z=0}) \ge  d .
\end{equation}
With the same arguments as in the proof of Proposition \ref{p:10}, this shows that $\Delta_{n-1}(f) \ge d$.

Let's now consider the case $(\alpha, \beta) = (0,1)$ with $\Delta_{n-1}(f) < 2d-2$. From the inequality (\ref{eq:le5}) at the end of the proof of  Proposition \ref{p:01} we deduce that 
$\Sigma^\ity_f$ is a reduced projective line with generic transversal type $A_1$. Note that, since 
the singular locus $\Sigma^\ity_f$ of $\{f_d=0\}$ is a line in $\bP^{n-1}$, one must have $n\ge 3$ and $d\ge 2$. This situation is treated in Theorem \ref{t:lineinfinity}. This result, proved independently in \S \ref{s:line}, tells that for such polynomials the top Betti defect is $\ge 2(n-1)(d-2) +1 > 2d-3$.

\subsection{Proof of Part (a)} \
 If  $f$ is general-at-infinity then 
$f$ has at most isolated singularities and $\Sigma_f^\ity = \emptyset$ by definition.  One applies formula (\ref{eq:F}) to get $b_{n-1}(f) = (d-1)^n$, since all the involved Milnor numbers are zero. This fact has already been observed in \cite{ST-gendefo}. 

Let us now prove the reciprocal with help of the above results. 
If the polynomial $f$ is not general-at-infinity then  $\Sigma^\ity_f \not= \emptyset$, and this implies in particular that $d = \deg f$ is greater than 1. If $\dim \Sing f \le 0$ and $\Sigma^\ity_f$ is a set of isolated points then we may use formula (\ref{eq:F}), in which there occurs at least one non-zero Milnor number, to conclude that the defect $\Delta_{n-1}(f) \ge 1$.
In case the singularities are non-isolated we use the above partial proof of part (b) of our theorem  to conclude that  $\Delta_{n-1}(f) \ge d-1 >0$.
\fin

%%%%%%%%%%%%%%%%%%%%%%%%%%%%%%%%%%
%%%%%%%%%%%%%%%%%%%%%%%%%%%%%%%%%%%%%%%%%%%%%%%%%%%%%%%%%%%%%%%%%%%%%
\section{Euler characteristic of projective hypersurfaces with one-dimensional singularities}\label{s:onedimatinfinity}
%Polynomials with one-dimensional singularities at the intersection with infinity
 
 We develop here a new method for computing the Euler characteristic of a projective hypersurface with non-isolated singularity.
We abute to a general result, Theorem \ref{t:euler_sum}.  We shall use its Corollary \ref{c:euler_sum} to build the proof of Theorem \ref{t:lineinfinity} in \S \ref{s:line}.

%\begin{lemma}\label{l:3}
 %There exists a deformation $f_s$ of $f$ with affine isolated singularities, such that $\Sigma^\ity_{f_s} = %\Sigma^\ity_f$ and that the hypersurface $\{ ({f_s})_{d-1} = 0\}$ is general position with respect to %$\Sigma^\ity_f$.
%\end{lemma}

%\begin{proof}
%Let us choose some homogeneous polynomial $h_{d-1}$ of degree $d-1$,  general enough such that the %nonsingular projective hypersurface $\{ h_{d-1}=0\}\subset \bP^{n-1}$ intersects $\Sigma^\ity_{f_s}$ %transversally at regular points.
%Taking the deformation $f_s = f+ s h_{d-1}$ we claim that the hypersurface $\{(f_s)_{d-1} = 0\}$ intersects %$\Sigma^\ity_{f_s}$ transversally at regular points. Indeed, 
%we have $\Sigma^\ity_{f_s} = \Sigma^\ity_f$ and the deformation $f_{d-1} + sh_{d-1}$ of the degree $d-1$ %part of $f$ is general for $s\not= 0$ in some neighbourhood of $0$. Then we are done since, by eventually %applying Lemma \ref{l:2} to $f_s$, which does not change the degrees $d$ and $d-1$ parts whenever $d>2$, we %may also assume that  $\Sing f_s$ has dimension zero. In case $d=2$, $h_{d-1}$ is linear and satisfies both %requirements, provided it is general enough.
%\end{proof}

Let $V := \{ f_d = 0 \}$ denote a hypersurface in $\bP^{n-1}= H^\ity$ of degree $d$ with singular locus $\hat \Sigma$ of dimension one, more precisely $\hat \Sigma$ consists of a union $\Sigma$ of irreducible curves and eventually a finite number of points $\{R_1,\ldots , R_\delta \}$.   Let $h_d$ be a general-at-infinity homogeneous polynomial of degree $d$ and consider the deformation $f_\e = f + \e h_d$. This is general-at-infinity for $\e \not= 0$ in some small enough disk centered at 0, by Proposition \ref{p:general}.
   For any $\e\in \bC$, let $V_\e := \{ f_{\e, d} := f_d + \e h_d= 0 \}$ be a pencil of projective hypersurfaces.

The genericity of $h_d$ ensures that $V_\e$ is nonsingular for all $\e \not=0$ in a small enough disk $\Delta \subset \bC$ centered at the origin. Let us consider the total space of the pencil:
\[  \bV_\Delta := \{ f_d + \e h_d = 0\} \subset \bP^{n-1} \times \Delta
\]
as germ at $V_0$
and the projection $\pi : \bV_\Delta \to \Delta$. We denote by $A = \{ f_d=h_d=0\}$ the axis of the pencil. One considers the polar locus of the map $(h_d, f_d) : \bC^n \to \bC^2$ and since this is a homogeneous set one takes its image in $\bP^{n-1}$ which will be denoted by $\Gamma(h_d, f_d)$. 

%%%%%%%%%%%%%
Let us spell out more precisely the meaning of ``general'' for $h_d$.
By using the Veronese embedding of degree $d$ we find a Zariski open set $\cO$ of linear functions in the target such that whenever $g\in \cO$ then its pull-back is a general homogeneous polynomial $h_d$ defining a hypersurface $H := \{ h_d = 0\}$ which is transversal to $V$ in the stratified sense, i.e. after endowing $V$ with some Whitney stratification, of which the strata are as follows:  the isolated singular points  $\{R_1\},\ldots , \{R_\delta\}$ of $V$ and the point-strata $\{Q_1\}, \ldots , \{Q_\gamma\} \subset \Sigma$, the components of $\Sigma \m \{ Q_1, \ldots , Q_\gamma \}$ and the open stratum $V\m \hat\Sigma$. Such $h_d$ will be called \textit{general}. This definition implies that $A$ intersects $\hat\Sigma$ at general points, in particular does not contain any point $Q_i$ or $R_i$. Moreover:
\begin{equation}  \label{eq:empty}
 \mbox{the germ of the polar locus } \Gamma_p(h_d, f_d)\subset \bP^{n-1} \mbox{ is empty at any point } p\in A\times \{ 0\}.
\end{equation}
 Let us explain and prove this claim. The germ $\Gamma_p(h_d, f_d)$ at some point $p$ is obtained as the polar locus of the localised map  $(\hat h_d, \hat f_d) : \bC^{n-1} \to \bC^2$ obtained by dehomogenisation.
At the regular points of $A= V\cap H$ the claim is due to the transversality of the hypersurfaces $V$ and $H$. At the points of $A\cap \Sigma$, we reason as follows: locally, the hypersurface $V$ is defined as the zero locus of the localised function $\hat f_d : \bC^{n-1} \to \bC$. One may suppose from the beginning that the fixed Whitney stratification of $V$ satisfies Thom condition\footnote{it has been proved that this is the case in general, see e.g. \cite{BMM}.} ($a_{\hat f_d}$). Since $h_d$ was chosen general, the hypersurface $H$ slices $\hat \Sigma$ only at the 1-dimensional Thom strata of $V$ (which are included in $\Sigma$) and the intersection is transversal. This implies that $H$ is transversal to the limits of the tangent spaces at the fibres of $f_d$, which means that the germ of the polar locus $\Gamma_p(\hat h_d, \hat f_d) = \Gamma_p(h_d, f_d) $ is empty, which was our claim. 

\begin{figure}[hbtp]
\begin{center}
%\epsfxsize=5cm
%\leavevmode
\resizebox{!}{2in}{\includegraphics{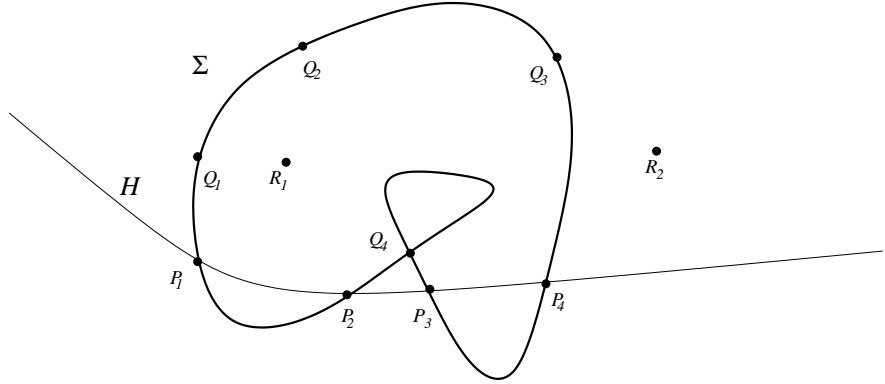}}
\end{center}
\caption{{\em Point-strata of $V$ and the intersection $H\cap \Sigma$
 }}
\label{f:2}
\end{figure}

With these preliminaries, we may prove the following:
\begin{lemma}\label{l:singlocus}
 The space $\bV_\Delta$ has isolated singularities: $\Sing \bV_\Delta = (A \cap \Sigma)
 \times \{ 0\}$, 
and  $\pi : \bV_\Delta \to \Delta$ is a map with 1-dimensional singular locus: $\Sing (\pi) = \hat \Sigma \times \{ 0\}$.
\end{lemma}
\begin{proof} 
The germ at $A$ of the singular locus
\[ 
 \Sing \bV_\Delta = \{ \partial f_d + \e \partial h_d = 0\} \cap A\times \Delta
\]
 is by definition the union of the finite set $A \cap\Sigma$ and the locus where the partial derivatives  $\partial f_d$ and $\partial h_d$ are dependent, which is some subgerm of the polar locus $\Gamma_A(h_d, f_d)$ at $A$. But, as shown above, this polar locus germ is empty for general $h_d$. 

Next, the map  $\pi: \bV_\Delta \to \Delta$ is a holomorphic map on a space with isolated singularities such that $\Sing \bV_\Delta \subset \Sigma \times \{ 0\}$. 
Then $\Sing \pi = \{ \partial f_d + \e \partial h_d =0\}$ may be read as the union of the singularities of the hypersurfaces $V_\e$ for $\e \in \Delta$. Now since  $h_d$ is general, the hypersurfaces $V_\e$ are nonsingular for  small enough $\e\not= 0$, like in the proof of Proposition \ref{p:general}. 
This shows that $\Sing \pi = \Sing V = \hat \Sigma$.
\end{proof}

We shall use (or continue to use) the following notations:  $A\cap \Sigma = \{ P_1, \ldots , P_\nu\}$,  $\Sigma^* := \Sigma \m (\{ P_i\}_{i=1}^\nu \cup (\{ Q_j\}_{j=1}^\gamma)$, $\cN :=$ small enough  tubular neighbourhood of $\Sigma^*$, and
$B_i, B_j, B_k$ are small enough Milnor balls within $\bV_\Delta \subset \bP^{n-1}\times \Delta$  at the points $P_i, Q_j, R_k$, respectively.

Since the Euler characteristic $\chi$ is a constructible functor, we have the following decomposition into a sum:
\begin{equation}\label{eq:chi}
 \chi(\bV_\Delta, V_\e) = \chi(\cN, \cN \cap V_\e) + \sum_{i=1}^\nu \chi(B_i, B_i\cap V_\e)+ \sum_{j=1}^\gamma \chi(B_j, B_j\cap V_\e)+ \sum_{k=1}^\delta \chi(B_k, B_k\cap V_\e)
\end{equation}

 The pair $(B, B\cap V_\e)$ in all of the above three sums represents the local Milnor data of a hypersurface germ of dimension $n-2$ in a space of dimension $n-1$. The last one $(B_k, B_k\cap V_\e)$ corresponds to the isolated hypersurface singularity of $V$ at $R_k$ with Milnor number $\mu_k \ge 1$, of which $\pi$ is a smoothing,  and therefore we have:
\[
 \chi(B_k, B_k\cap V_\e) = (-1)^{n-1} \mu_k
\]

%We shall prove that all the nonzero terms contribute with the same sign in the formula %(\ref{eq:chi}).

For the first term, since the map $\pi : \bV_\Delta \to \Delta$  has a trivial transversal  structure along $\Sigma_r^*$, where $\Sigma_r$ is some irreducible component in the decomposition of $\Sigma$, we have the equality:
\[
 \chi(\cN, \cN \cap V_\e) = \sum_r \chi(\Sigma_r^*) \chi(B_r, F_r^\pitchfork).
\]
where $(B_r, F_r^\pitchfork)$ is the transversal Milnor data at some point of $\Sigma_r^*$, namely $B_r$ is a Milnor ball of the transversal singularity and $F_r^\pitchfork$ is the corresponding transversal Milnor fibre. Note that this is the Milnor data of an isolated hypersurface singularity of dimension $n-3$; its  Milnor number will be denoted by $\mu_r^\pitchfork$ and this does not depend on the choice of the point on $\Sigma_r^*$. We therefore have: $\chi(B_r, F_r^\pitchfork)= (-1)^{n-2}\mu_r^\pitchfork$.

 We also have $\chi(\Sigma_r^*) = 2-2g_r - \nu_r -\gamma_r$ where $g_r$ is the genus of $\Sigma_r$, and where
$\nu_r$ and $\gamma_r$ are the numbers of points $P_i$ and $Q_j$ on $\Sigma_r$, respectively. Then:
\begin{equation}\label{eq:chiN}
 \chi(\cN, \cN \cap V_\e)= (-1)^{n-1} \sum_r (\nu_r+ \gamma_r + 2g_r -2)\mu_r^\pitchfork .
\end{equation}

Let us show that the contribution of the axis in the formula (\ref{eq:chi}) is null.

\begin{lemma}\label{l:B_i}
 $\chi(B_i, B_i\cap V_\e) = 0$.
\end{lemma}
\begin{proof}
Let $p\in A \cap \Sigma$ be the center of the Milnor ball $B_i$. Consider the map $(\pi, h_d) :  B_i \to \Delta \times \Delta'$. Consider the germ of the polar locus of this map at $p$, denoted by $\Gamma(\pi, h_d)$.
Then it follows from the definition of the polar locus that some point $(x,\e) \in \bV_\Delta$, where $\e = - f_d(x)/h_d(x)$, is contained in  $\Gamma(\pi, h_d) \setminus (\{ f_d=0\} \cup \{ h_d=0\})$ if and only if $x \in \Gamma(f_d, h_d) \setminus (\{ f_d=0\} \cup \{ h_d=0\})$.
But we have shown in (\ref{eq:empty}) that $\Gamma(f_d, h_d)$ is empty at $p$.
The absence of the polar locus implies that $B_i\cap V_\e$ is homotopy equivalent (by deformation retraction) to the space $B_i\cap V_\e \cap \{ h_d=0\}$. The latter is the slice by $\e =$ constant of the space $\bV \cap \{ h_d=0\} = \{ f_d=0\} \times \Delta$, which is a product space. Since this is homeomorphic to the complex link of this space and a product space has contractible complex link, we deduce that $B_i\cap V_\e$ is contractible too. Since $B_i$ is contractible itself, we get our claim.
\end{proof}

%%%%%%%%%%%%%%%%%%%%%%%%%%
%%%%%%%%%%%%%%%%%%%%%%%%%%

Let $V_\e := \{ f_d + \e h_d =0\}$ be the general pencil considered above, $\e \ne 0 \in \Delta$.  
We have  $ \chi (\bV_\Delta, V_\e) = \chi(V) - \chi^{n-1,d}$  since $V_\e$ is a general hypersurface of degree $d$ in $\bP^{n-1}$ and since $\bV_\Delta$ retracts to its central fibre $V = \{ f_d = 0 \}$.  
Then the preceding considerations prove the following:
\begin{theorem}\label{t:euler_sum} 
Let $V := \{ f_d = 0 \}\subset \bP^{n-1}$ be a hypersurface of degree $d$ where $\Sing V$ is a union of curves and isolated points.   
Then, in the above notations:
 \begin{equation}\label{eq:chi_sum0}
   \chi(V) = \chi^{n-1,d} + (-1)^{n-1} \sum_r (\nu_r+ \gamma_r + 2g_r -2)\mu_r^\pitchfork +
 \sum_{j=1}^\gamma \chi(B_j, B_j\cap V_\e)+ (-1)^{n-1}\sum_{k=1}^\delta \mu_k.
\end{equation}
\fin
\end{theorem}

Our aim is to compute, via formula (\ref{eq:B}), the top Betti defect of polynomials $f$ with  $\dim \Sing f \le 0$ and  $\Sigma^\ity_f$ is a union of curves and isolated points,  such that 
$\dim \Sigma^\ity_f \cap \{ f_{d-1} = 0\} \le 0$.
We are considering the deformation $f_\e = f + \e h_d$, where $h_d$ is a general-at-infinity homogeneous polynomial of degree $d$. By Proposition \ref{p:general},  $f_\e$ is a general-at-infinity polynomial for $\e \not= 0$. In the notations of Proposition \ref{p:defect},  $\Delta\chi^{\infty}= \chi (\bV_\Delta, V_\e)$.
Then Theorem \ref{t:euler_sum} reads: 
 
\begin{corollary}\label{c:euler_sum} 
Let $f:\bC^n\to \bC$ be a polynomial of degree $d$ with $\dim \Sing f \le 0$ and  $\Sigma^\ity_f$ is a union of curves and isolated points. 
Then:
\begin{equation}\label{eq:chi_sum}
    \Delta\chi^{\infty} = (-1)^{n} \sum_r (\nu_r+ \gamma_r + 2g_r -2)\mu_r^\pitchfork -
 \sum_{j=1}^\gamma \chi(B_j, B_j\cap V_\e)+ (-1)^{n}\sum_{k=1}^\delta \mu_k
\end{equation}

\fin
\end{corollary}

The only part of the formula (\ref{eq:chi_sum}) which is not explicitly computed is the sum of  $\chi(B_j, B_j\cap V_\e)$ which runs over the Whitney point-strata $Q_j$ of the hypersurface $V$. One may compute it in  particular cases, as we show in the following section.

%%%%%%%%%%%%%%%%%%%%%%%%%%%%%%%%%%%%%%%%%
%%%%%%%%%%%%%%%%%%%%%%%%%%%%%%%%%%%%%%%%%%

\section{Polynomials with line singularities at infinity}\label{s:line}
 
A natural class of polynomials is the one where $\Sigma^\ity_f$ is a reduced line with Morse generic transversal type. Indeed, by the proof of Proposition \ref{p:01}, if the top Betti defect of $f$ is between $d$ and $2d-3$, then
$f$ might have such type of singularities. We prove here Theorem \ref{t:lineinfinity}, namely that
if $f$ has at most isolated affine singularities and $\Sigma^\ity_f$ is a reduced projective line with Morse generic transversal then $\Delta_{n-1}(f) \ge 2(n-1)(d-2) +1$. 

As remarked in \S \ref{ss:end1}, the existence of such singularities implies $n\ge 3$. For $n=3$ our formula reads: $\Delta_{n-1}(f) \ge 4d -7$ which shows that our result specialises to the estimation proved in \cite{ALM} for a particular class of polynomials in  $3$ variables, with $\dim\Sing f \le 0$ and $\dim \Sigma^\ity_f =1$, and with no singularities at infinity in the sense of \cite{ST}.

%%%%%%%%%%%%%%%%%%%%%%%%%%%%%%%
\subsection*{Proof of Theorem \ref{t:lineinfinity}}
By eventually deforming the  $d-1$ homogeneous part of $f$ we get that the intersection  $\Sigma \cap  \{ f_{d-1}=0\}$ is of dimension $\le 0$.
Using the definitions in \S \ref{s:onedimatinfinity} in our particular setting, we have by assumption  $\delta =0$, $r=1$, $g_1=0$, $\mu_1^\pitchfork =1$,   and $\nu_1 = \nu = \mult(\Sigma, \{ h_d = 0\}) = \deg h_d =d$.   We apply formula (\ref{eq:B}) and Corollary \ref{c:euler_sum} and we get:
\begin{equation}\label{eq:defectfinal}
\Delta_{n-1}(f) = \sum_{p\in\Sigma_f^\ity \cap \{ f_{d-1}=0\}}\mu_p(\overline{X_0}) +  (d + \gamma -2) +  (-1)^{n-1}\sum_{j=1}^\gamma \chi(B_j, B_j\cap V_\e).
\end{equation}
Let us  first evaluate the sum of Milnor numbers of the general fibre $X_0$ of $f$. For a general $d-1$ homogeneous part of $f$ we get that the intersection  $\Sigma \cap  \{ f_{d-1}=0\}$ consists of $d-1$ simple points, each of which being an $A_1$ singularity of $X_0$.
 This implies that the above first sum is bounded from below by $d-1$.

The number $\gamma$ counts the special points $Q_j$ on the singular line $\Sigma :=  \Sigma_f^\ity$. Then the last sum has $\gamma$ terms and we need to determine for each of them the contribution $\chi(B_j, B_j\cap V_\e)$. For this we need the deformation theory of line singularities, founded by Siersma \cite{Si1} and subsequently developed by several authors. Let us assume without loss of generality that the line $\Sigma \subset \bP^{n-1}$ is the zero locus of the ideal $I= (x_1, \ldots, x_{n-2})$. We remark first that the ideal of homogeneous polynomials $g : (\bC^{n},0) \to (\bC,0)$ such that $\Sing g \supset \Sigma$ is spanned by the polynomials of the form $g(x) = \sum_{i, j=1}^ {n-2} h_{ij}(x) x_ix_j$ where  $h_{ij}(x)$ are polynomials depending on all variables $x_1, \ldots , x_{n}$.
This was established in \cite{Si1} and \cite{Pe1} in the germ case; our situation is slightly different but the same proof applies.

%Our present situation is slightly different but the same proof applies. Let $g := \hat f_d : %\bC^{n-1} \to \bC$ be the localisation of the polynomial $f_d$ in the chart $x_n =1$. Instead of a %germ, we consider here the global fibre $g^{-1}(0)$. Then the primitive ideal $\int_I$ is equal to %$I^2$ and $g\in \int_I$ is equivalent to $\Sigma \subset \Sing g^{-1}(0)$.

  In our setting the functions $h_{ij}(x)$ are of degree $\le d-2$.  
Following the deformation theory in  \cite{Si1}, by deforming $h_{ij}(x)$
we get a generic \textit{transversal Hessian} $\cH(x) := det (h_{ij})_{i,j}(x)$ and this implies that $g$ has generic singularity type $A_\ity$ along $\Sigma$. Then the point-strata $Q_j$ are precisely the type $D_\ity$ singularities. Following Siersma's theory \cite{Si1}, the number of $D_\ity$ points is equal to the degree of the Hessian $\cH(0, \ldots, 0, x_n)$. In the generic case this degree turns out to be equal to $(d-2)(n-2)$. We may then take this value of $\gamma$ in the formula (\ref{eq:defectfinal}) as a minimum. 
We also get from \cite{Si1} that the Milnor fibre of a  $D_\ity$-singularity is homotopy equivalent to a $(n-2)$-sphere, therefore $\chi(B_j, B_j\cap V_\e) = (-1)^{n-1}$.

Finally, putting together the lower bounds we get:
\[ 
 \Delta_{n-1}(f) \ge  d-1 + d + 2(d-2)(n-2) - 2 = 2(n-1)(d-2) +1.
\]
\fin

\begin{remark}\label{r:ex}
  The proof actually shows that whenever $h_{ij}$ and $f_{d-1}$ are generic we have the equality $\Delta_{n-1}(f)  = 2(n-1)(d-2) +1$, which means that the bound of Theorem \ref{t:lineinfinity} is sharp for any $n \ge 3$. An explicit example for $n=3$ is given in \S \ref{s:intro}, Example \ref{ex4}.
\end{remark}

%%%%%%%%%%%%%%%
%

%%%%%%%%%%%%%%%%%%%%%%%%%%%%%%%%%%%%%%%%%%%%%%%%%%%%%%%%%%%%%%%%%%%%%

\end{document}